\documentclass[letterpaper, 10 pt, conference]{ieeeconf}
\IEEEoverridecommandlockouts

\usepackage{amsmath,amssymb,amsfonts}
\usepackage{graphicx}
\usepackage{textcomp}
\usepackage{preamble}
\usepackage{xcolor}
\usepackage{pgfplots}
\usepackage{pgfplotstable}
\usepackage{float}
\usepackage{wrapfig}
\usepackage{caption}
\usepackage{subcaption}
\renewcommand{\textfloatsep}{0.15 cm}%

\pgfplotsset{compat=1.7}

\definecolor{linkblue}{RGB}{11,104,149}
\hypersetup{
    colorlinks=true,
    linkcolor=linkblue,
    urlcolor=linkblue,
}

\pgfplotsset{
SmallBarPlot/.style={
    font=\small,
    ybar,
    width=\linewidth,
    ymin=0,
    xtick=data,
    bar width=0.4em,
},
}

\usepackage[style=ieee, url=false, doi=false, isbn=false, eprint=false]{biblatex}
\begin{document}
\title{PANTR: A proximal algorithm with trust-region updates for nonconvex constrained optimization}
\author{Alexander Bodard,\; Pieter Pas\; and\; Panagiotis Patrinos
\thanks{
    This work was supported by:
    Research Foundation Flanders (FWO) PhD grant No. 11M9523N and research projects G081222N, G033822N, G0A0920N;
    Research Council KU Leuven C1 project No. C14/18/068;
    EU's Horizon 2020 research and innovation programme under the Marie Skłodowska-Curie grant agreement No. 953348;
    Fonds de la Recherche Scientifique - FNRS and the Fonds Wetenschappelijk Onderzoek - Vlaanderen under EOS project no 30468160 (SeLMA).\newline
    \indent 
    The authors are with the Department of Electrical Engineering (ESAT-STADIUS), 
    KU Leuven, Kasteelpark Arenberg 10, 3001 Leuven, Belgium.
    Email: \texttt{\fontsize{7}{7}\selectfont\{alexander.bodard,pieter.pas,panos.patrinos\}@esat.kuleuven.be}
}
}

\maketitle
\thispagestyle{empty} %

\begin{abstract}
This work presents \pantr{}, an efficient solver for nonconvex constrained optimization problems, that is well-suited as an inner solver for an augmented Lagrangian method.
The proposed scheme combines forward-backward iterations with solutions to trust-region subproblems:
the former ensures global convergence, whereas the latter enables fast update directions.
We discuss how the algorithm is able to exploit exact Hessian information of the smooth objective term through a linear Newton approximation, while benefiting from the structure of box-constraints or $\ell_1$-regularization.
An open-source \Cpp{} implementation of \pantr{} is made available as part of the NLP solver library \alpaqa{}.
Finally, the effectiveness of the proposed method is demonstrated in nonlinear model predictive control applications.
\end{abstract}

\section{Introduction}
\label{sec:introduction}

\subsection{Background and motivation}

Various areas of science and engineering naturally give rise to constrained, potentially nonconvex optimization problems, motivating the need for efficient solvers.
A prominent example is found in the field of model predictive control (MPC) \cite{rawlings_model_2017}.
Historically, there has been a strong focus on linear MPC, where linear system dynamics lead to convex (and often linear-quadratic) optimization problems. However, recent advances in hardware have sparked an increased interest in nonlinear MPC (NMPC), which produces more challenging nonconvex, nonlinear programs (NLPs).

In this context, not only the efficiency of the solver, but also its reliability and memory requirements are critical.
Classical techniques for solving the NLPs include interior point (IP) methods and sequential quadratic programming (SQP) \cite{nocedal_numerical_2006}.
State-of-the-art IP solvers such as \ipopt{} \cite{wachter_implementation_2006} are reliable as general-purpose solvers, but have the disadvantage that they do not easily exploit warm starts and have large memory requirements.
SQP involves repeatedly solving a quadratic program (QP), and relies on efficient QP solvers.
An overview of recent advances in this area is given in \cite{kouzoupis_recent_2018}.
As an alternative to these classical approaches, \panoc{} \cite{stella_simple_2017} is a first-order method that combines forward-backward (FB) iterations with quasi-Newton directions to attain fast local convergence.
\Panoc{} benefits from warm starts and has a much smaller memory footprint.
Since \panoc{} was originally proposed, its effectiveness in real-time MPC applications has been shown in various works \cite{sathya_embedded_2018-1,small_aerial_2019,lindqvist_reactive_2022}.

The underlying premise of \panoc{} is that the proximal mapping of the objective's nonsmooth term can be efficiently evaluated (\S \ref{sec:problem-statement}).
To handle more general constraints, the \alpaqa{} package \cite{pas_alpaqa_2022} implements an augmented Lagrangian method (ALM) that uses \panoc{} as an inner solver. 
Both quasi-Newton and Gauss-Newton \cite{pas_gauss-newton_2022} variants of \panoc{} speed up convergence through approximate second-order information.
In contrast, the present work applies directions generated by a linear Newton approximation (LNA) \cite{facchinei_finite-dimensional_2003}, with the goal of enabling faster convergence than \panoc{} through \emph{exact Hessian information} of the smooth term of the cost.
For nonconvex optimization problems, the LNA may have negative eigenvalues.
We therefore propose a semismooth Newton scheme that computes update directions as solutions to trust-region (TR) subproblems, allowing the exploitation of directions of negative curvature of the LNA.

\subsection{Contributions}

The contributions of this work include:
(\emph{i}) we propose \pantr, a novel proximal method that combines FB iterations with solutions to TR subproblems and as such is able to exploit exact second-order information for nonconvex problems;
(\emph{ii}) we prove global convergence of \pantr, without requiring feasibility of the iterates;
(\emph{iii}) we show that if the TR constraint becomes inactive, \pantr{} generates fast Newton directions;
(\emph{iv}) an open-source \Cpp{} implementation of \pantr{} is made available as part of \alpaqa{};\footnote{\url{https://github.com/kul-optec/alpaqa}}
(\emph{v}) we show the effectiveness of \pantr{} on a number of NMPC problems. In particular, \pantr{} appears to perform well on problems where \panoc{} struggles.

\subsection*{Notation}

By $\N$, $\R$ and $\Rbar = \R \cup \{ +\infty \}$ we denote the set of natural, real and extended real numbers respectively.
The restriction of $\N$ to $[i, j]$ is written as $\N_{[i, j]} \eqdef \N \cap [i, j]$ and 
$x_i$ denotes the \(i\)'th component of $x \in \R^n$.
We write $x_\mathcal{I} = (x_i)_{i \in \mathcal{I}}$ for an index set $\mathcal{I} \subseteq \N_{[1, n]}$, and
denote the Euclidean inner product and norm by $\langle \cdot, \cdot \rangle$ and $\|\cdot\|$ respectively.
The proximal operator of a function $h : \R^n \to \Rbar$ is $\prox_{h} (x) \eqdef \argmin_u \{ h(u) + \frac{1}{2} \Vert u - x \Vert^2 \}$.
By $\Pi_C$ we denote the Euclidean projection on a set $C$.
The set $\fix{T} \eqdef \{ x \in \R^n \mid x \in T(x) \}$ contains the fixed points of an operator $T : \R^n \rightrightarrows \R^n$.
We denote the set of $k$ times continuously differentiable functions by $\mathcal{C}^k$.
We say that $f \in \mathcal{C}^1$ is $L_f$-smooth when $\nabla f$ is $L_f$-Lipschitz continuous.
We denote the Clarke generalized Jacobian of a function $F : \R^n \to \R^m$ by $\partial_C F$ \cite{clarke_optimization_1990}.
$lsc$ and $osc$ refer to lower and outer semicontinuity as in \cite{rockafellar_variational_1998}.

\section{Problem statement and preliminaries} \label{sec:problem-statement}

\Alpaqa{} \cite{pas_alpaqa_2022} combines an augmented Lagrangian method with inner solvers like \panoc{} \cite{stella_simple_2017,de_marchi_proximal_2022} to tackle NLPs %
\begin{equation} \tag{P} \label{eq:original-problem}
    \begin{aligned}
        &\minimize_{x \in \R^n} &&f(x)\\
        &\stt && \underline{x} \leq x \leq \overline{x}, &&& \underline{z} \leq g(x) \leq \overline{z},
    \end{aligned}
\end{equation}
where $f : \R^n \to \R$ and $g : \R^n \to \R^m$ may be nonconvex.
\Panoc{} attains fast convergence in solving the inner problem \eqref{eq:problem-statement} through quasi-Newton or Gauss-Newton directions \cite{pas_gauss-newton_2022}.
However, disadvantages of the existing methods using quasi-Newton and Gauss-Newton directions are that 
(\emph{i}) the Hessian is kept positive definite, creating a potential disagreement with the curvature of the underlying function; 
(\emph{ii}) no \emph{exact} Newton directions are exploited.
Hence, we propose an alternative solver that exploits \emph{second-order structure} of the inner problem.

For nonconvex problems, using the exact Hessian or a possibly indefinite approximation thereof requires \textit{regularization} whenever the quadratic model is unbounded below. 
Techniques like Levenberg-Marquardt or modification of the eigenvalues require careful selection of the regularization parameter to attain fast convergence.
In this work, we choose to implicitly regularize by limiting the norm of the Newton step, which corresponds to solving a TR subproblem.

The remainder of this section briefly reviews the augmented Lagrangian method (ALM), existing FB schemes like \panoc{}, and TR methods.

\subsection{Augmented Lagrangian method} \label{sec:alm}

Define $C \eqdef \{ x \mid \underline{x} \leq x \leq \overline{x} \}$ and $D \eqdef \{ z \mid \underline{z} \leq z \leq \overline{z} \}$.
By use of a slack vector $z$, Problem \eqref{eq:original-problem} can be written as
\begin{align} \tag{P-ALM} \label{eq:alm}
        &\minimize_{x \in C, z \in D} &&f(x) &&&\stt &&&&z = g(x).
\end{align}
Given a positive definite diagonal matrix $\Sigma \in \R^{m \times m}$, we define the augmented Lagrangian function with penalty factor $\Sigma$ as $\lagr_\Sigma(x, z, y) \eqdef f(x) + \langle y, g(x) - z \rangle + \frac{1}{2} \Vert g(x) - z \Vert_\Sigma^2$.
ALM applied to Problem \eqref{eq:alm} iteratively 
(\emph{i}) minimizes $\lagr_\Sigma$ w.r.t.\,$x$ and $z$; 
(\emph{ii}) updates the Lagrange multipliers $y$; and 
(\emph{iii}) updates the penalty factors $\Sigma$.
For more detailed information on ALM in the context of \alpaqa, see \cite{pas_alpaqa_2022}.

\subsection{Forward-backward schemes}

Consider the composite minimization
\begin{equation} \label{eq:problem-statement} \tag{P-FB}
    \minimize_{x \in \R^n} \quad\varphi(x) \eqdef \f(x) + \g(x)
\end{equation}
of two functions ${\f : \R^n \to \R}$ and ${\g : \R^n \to \Rbar}$ where $\prox_{\gamma \g}$ is efficiently evaluated.
Throughout this work, we assume that (\emph{i}) $\f \in \mathcal{C}^1$ is $\lipschf$-smooth; (\emph{ii}) $\g$ is proper, lsc and convex; and (\emph{iii}) $\varphi \equiv \f + \g$ is lower bounded.
These are standard assumptions for first-order splitting schemes, see e.g. \cite[Assumption 10.1]{beck_first-order_2017}.\footnote{
    Section \ref{sec:pantr-adaptive} relaxes the first assumption to include all functions $\f$ with \emph{locally} Lipschitz-continuous gradients.
}
Note that \eqref{eq:problem-statement} reduces to \emph{nonconvex constrained optimization} by defining $\g$ as the indicator function of a constraint set.
Also observe that the ALM inner problem, i.e.\,the minimization of $\lagr_\Sigma$ with respect to $x$ and $z$, can be formulated in the form \eqref{eq:problem-statement} by defining $\f(x) = f(x) + \frac{1}{2} \mathbf{dist}_\Sigma^2 (g(x) + \Sigma^{-1} y, D)$ and $\g = \delta_C$, where $\mathbf{dist}_\Sigma (\cdot, D) \eqdef \minbold_{z \in D} \left\{ \Vert z - \cdot \Vert_\Sigma \right\}$.\footnote{
    Since $\mathbf{dist}_\Sigma^2(\cdot, D)$ has a piecewise linear gradient, it is Lipschitz-smooth.
    Hence, for \emph{locally} Lipschitz-smooth functions $g$, also $\f$ is \emph{locally} Lipschitz-smooth.
}

Solutions to \eqref{eq:problem-statement} are fixed points of the \emph{forward-backward operator}  
\(T_\gamma(x) \eqdef \prox_{\gamma \g}(x - \gamma \nabla \f(x))\),
or equivalently, zeros of the \emph{fixed-point residual} $R_\gamma(x) \eqdef \nicefrac{1}{\gamma} \left( x - T_\gamma(x) \right)$, where \( \gamma > 0 \).
The forward-backward splitting (FBS) scheme iteratively applies $T_\gamma$, and converges to a fixed-point of $T_\gamma$ for sufficiently small step sizes \(\gamma\).
Various methods have been proposed to accelerate the convergence of FBS.
For example, both \panoc{} \cite{stella_simple_2017,de_marchi_proximal_2022} and \zerofpr{} \cite{themelis_forward-backward_2018} aim to find zeros of $R_\gamma$ through a quasi-Newton line search procedure that uses the forward-backward envelope (FBE) \cite{patrinos_proximal_2013,stella_forwardbackward_2017,themelis_forward-backward_2018} as a merit function.
\begin{definition}[FBE]\label{def:FBE}
    The FBE of $\varphi$ with parameter $\gamma > 0$ is
    \begin{equation*}
    \varphi_\gamma(x)=\infbold_{u \in \R^n}\f(x) + \langle \nabla \f(x), u - x \rangle + \g(u)+\tfrac{1}{2 \gamma}\|u-x\|^2.
    \end{equation*}
\end{definition}
\noindent
The FBE is a useful metric, since its minimization is equivalent to the minimization of $\varphi$ in the sense that $\argmin \varphi \equiv \argmin \fbe$ for $\gamma \in \left( 0, \nicefrac{1}{\lipschf} \right)$ \cite[Pr. 2.3 iii]{stella_forwardbackward_2017}.
Remark that when $\f \in \mathcal{C}^2$, the FBE is continuously differentiable with $\nabla \fbe(x) = Q_\gamma(x) R_\gamma(x)$ where $Q_\gamma(x)\eqdef\id-\gamma \nabla^2 \f(x)$ \cite[Th. 2.6]{stella_forwardbackward_2017}.
Yet, $\fbe$ is not twice differentiable in general, since $R_\gamma$ is not generally differentiable.

\subsection{Trust-region methods}\label{sec:tr}

To minimize a function $F : \R^n \to \R$, TR methods define a model $m_k$ that locally approximates the objective function $F$, typically based on a second order Taylor expansion. 
A TR step $d_k$ is then computed by minimizing the model $m_k$ over all points within a distance $\trradius_k$ of the current iterate, i.e.\,
\begin{equation} \label{eq:tr-general}
    \hspace{-1.1em}
    \begin{aligned}[b]
        &\minimize_{d} &&m_k(d) \eqdef F(x_k) + \langle \nabla F(x_k), d \rangle + \tfrac{1}{2} \langle B_k d, d \rangle\\
        &\stt &&\Vert d \Vert \leq \trradius_k
    \end{aligned}\hspace{-2em}
\end{equation}
where $B_k$ is the (approximate) Hessian matrix of $F$ evaluated at $x_k$.
The TR radius $\trradius_k$ is updated depending on how well $m_k$ approximates $F$, through the ratio
\begin{equation}\label{eq:classical-tr-ratio}
    \trratio_k = \tfrac{F(x_k) - F(x_k + d_k)}{m_k(0) - m_k(d_k)}.
\end{equation}
The TR subproblem can be interpreted as adaptively regularizing $B_k$, since for any global minimizer of \eqref{eq:tr-general} \cite[Cor. 7.2.2]{conn_trust_2000}
\begin{equation} \label{eq:tr-regularization}
    \exists \lambda \geq 0: (B_k\!+\!\lambda \id)\,d = -\nabla F(x_k)\;\;\text{and}\;\; B_k\!+\!\lambda \id \succeq 0.
\end{equation}
Approximately solving this subproblem efficiently is critical for the overall performance of TR methods.
The \emph{Steihaug conjugate gradient (CG) method} \cite{steihaug_conjugate_2006} is widely used, although other techniques exist.
The interested reader is referred to \cite{conn_trust_2000} for an extensive study of TR methods.

\section{\pantr} \label{sec:pantr}

This section introduces the \pantr{} scheme for solving \eqref{eq:problem-statement}. First, the semismooth Newton system for the root-finding problem $R_\gamma(x^\star) = 0$ is transformed into a minimization problem, and it is regularized by fitting it into a TR framework. Next, we globalize the scheme and analyze its convergence properties.
Finally, we propose an adaptive step size procedure, such that no explicit knowledge of the Lipschitz constant $\lipschf$ is needed.

\subsection{Newtonian directions for the fixed-point residual} \label{sec:pantr-newton}

\newcommand\xnewt{{\widehat{x}}}
\newcommand\Rxnewt{R_\gamma(\xnewt)}
\newcommand\lnaRxnewt{\lna{R_\gamma}}

As a surrogate for solving \eqref{eq:problem-statement}, we are interested in finding a point $x_\star$ such that $R_\gamma(x_\star) = 0$.
\Panoc{} addresses this problem using a Newton-type approach and computes update directions $d_k = -B_k^{-1} R_\gamma(x_k)$, where $B_k$ is an invertible matrix.
If $B_k$ captures first-order information of $R_\gamma$, this enables superlinear convergence when sufficiently close to a strong local minimum.
The operators $B_k$ are typically defined using an L-BFGS scheme.
However, this enforces symmetry and positive definiteness of $B_k$, which restricts the algorithm's ability to exploit directions of negative curvature.

Ideally, we are interested in performing Newton's method, which defines $B_k$ as the exact Jacobian of the fixed-point residual $\jac{R_\gamma}(\xnewt)$.
Unfortunately, $R_\gamma$ is \emph{not in general differentiable}.
Hence, we propose to use a linear Newton approximation (LNA) \cite{facchinei_finite-dimensional_2003} of \(R_\gamma\) instead.
Define the generalized Hessian of $\f$ as $\partial^2 \f \eqdef \partial_C(\nabla \f)$. 
Then
\begin{equation} \label{eq:lna-R}
    \hat{\partial} R_\gamma(x) \eqdef \bigg\{ \tfrac{1}{\gamma} \left( \id - \mathcal{P_\gamma} \mathcal{Q_\gamma} \right) \bigg| \begin{subarray}
        \mathcal P_\gamma \in \,\partial_C \prox_{\gamma \g}(x - \gamma \nabla \f(x))\\ 
        \mathcal Q_\gamma \in \,\id - \gamma \partial^2 \f(x)
    \end{subarray}\bigg\}
\end{equation}
is a LNA of $R_\gamma$ at $x_\star$ under a mild semismoothness assumption \cite[Prop. 1]{pas_gauss-newton_2022}. 
Remark that this is a generalization of the LNA originally proposed in \cite[Prop. 4.13]{themelis_acceleration_2019} for $\f \in \mathcal{C}^2$.

\newcommand\Pg{\mathcal P_{\microspace\gamma}}
At a given iterate \(\xnewt\), using the residual \(\Rxnewt\) and its LNA $\lnaRxnewt \in \hat{\partial}R_\gamma(\xnewt)$, we thus aim to compute an update direction $d$ that solves the Newton system $\lna{R_\gamma} d = -\Rxnewt$.
To this end, the following theorem introduces a constrained minimization problem, the 
solution of which matches the solution to the Newton system of the fixed-point
residual under the idempotence assumption \(\Pg^2 = \Pg\).
Observe that this condition can be satisfied when $\g = \delta_S$ where $S$ is a polyhedral set \cite[\S 6.2b]{themelis_acceleration_2019} and, by the Moreau decomposition, when $\g$ is the support function of a polyhedral set $S$.
\begin{theorem}\label{thm:minimization-newton-equiv}
    {\normalfont(Solutions of Newton system)}
    Select \(H \in \partial^2\f(\xnewt)\) and \(\Pg \in \partial_C \prox_{\gamma \g}(\xnewt - \gamma \nabla \f(\xnewt))\),
    such that \(\lna{R_\gamma} = \tfrac{1}{\gamma} \left( \id - \Pg \mathcal{Q_\gamma} \right)\).
    Suppose that \(\Pg^2 = \Pg\), and define \(\Pg^\perp \eqdef I - \Pg\).
    Then, a solution \(d^\star\) of the minimization problem
    \begin{equation}\label{eq:minimize-eq-constr}\normalfont
        \begin{aligned}
            &\minimize_{d} && \tfrac12 \langle d, Hd \rangle + \langle \Rxnewt, d \rangle \\
            &\stt &&\Pg^\perp d = -\gamma \Pg^\perp \Rxnewt
        \end{aligned}\hspace{-2em}
    \end{equation}
    also solves the Newton system of the fixed-point residual,
    \begin{equation}\label{eq:newton-lna-system}
        \lnaRxnewt\, d^\star = -\Rxnewt.
    \end{equation}
    \begin{proof}
        Being the generalized Jacobian of a proximal mapping, \(\Pg\) is symmetric \cite[Thm.\,15.4.12]{themelis_acceleration_2019}.
        The stationarity condition of problem \eqref{eq:minimize-eq-constr} is given by
        \(Hd^\star + \Rxnewt + \Pg^\perp \lambda = 0\), for some multiplier \(\lambda\).
        Multiplying by \(\Pg\) and using \(\Pg\Pg^\perp=0\)
        yields \(\Pg Hd^\star = -\Pg \Rxnewt\). Replacing %
        \(\Pg=I-\Pg^\perp\) on the right-hand side
        and substituting the constraint results in
        \(\Pg Hd^\star = -\Rxnewt - \gamma^{-1} \Pg^\perp d^\star\),
        implying \eqref{eq:newton-lna-system}.
    \end{proof}
\end{theorem}

\subsection{Box constraints and \(\ell_1\)-norm}
\newcommand\J{\mathcal J}
\newcommand\K{\mathcal K}
\newcommand\Perm{\mathrm P}
\newcommand\Pkj{\Perm_{\!\K\J}}
\newcommand\Pkjt{\Pkj^\top}
\newcommand\dK{d_{\microspace\constact}}
\newcommand\dJ{d_{\microspace\constinact}}

We now specialize to two classes of functions whose proximal operators satisfy the idempotence assumption \(\Pg^2 = \Pg\), i.e.\,
the indicators of rectangular boxes, and the \(\ell_1\)-norm.
First, consider the case where $\g$ is the indicator of the set of box
constraints $C = \left\{ x \mid \underline{x} \leq x \leq \overline{x} \right\}$.
The proximal operator then reduces to a projection onto $C$, i.e.\, $\prox_{\gamma \g} = \Pi_C$ \cite[Sec. 1.G]{rockafellar_variational_1998}.
Similarly to the derivation in \cite[Sec.\,III]{pas_alpaqa_2022},
define the set $\K \eqdef \underline{\K}(\xnewt) \cup \overline{\K}(\xnewt)$ of indices -- for brevity we omit the argument of $\K$ -- corresponding to \emph{active} box constraints on $\xnewt - \gamma \nabla \f(\xnewt)$ (i.e.\,after a forward step  on $\xnewt$) where
\begin{equation}
    \begin{aligned}
        \underline{\K}(\xnewt) &\eqdef \left\{ i \in \N_{[1, n]} \mid (\xnewt)_i - \gamma \nabla_{x_i} \f(\xnewt) \leq \underline{x}_i \right\}, \\
        \overline{\K}(\xnewt) &\eqdef \left\{ i \in \N_{[1, n]} \mid \overline{x}_i \leq (\xnewt)_i - \gamma \nabla_{x_i} \f(\xnewt) \right\}.
    \end{aligned}
\end{equation}
We denote the complement of $\K$, containing the indices of the \emph{inactive} constraints, by $\mathcal{J} = \N_{[1, n]} \backslash \K$.
To simplify notation, consider a row permutation matrix $\Pkj$ that reorders the rows with indices $k \in \K$ before those with indices $j \in \mathcal{J}$.
In this case, we can select $\mathcal{P_\gamma} \in \partial_C(\prox_{\gamma \g}(\xnewt - \gamma \nabla \f(\xnewt)))$ to be $\mathcal{P_\gamma} = \Pkjt \left(\begin{smallmatrix}
    0_{|\microspace\constact\microspace|} & 0_{\phantom{|\microspace\constinact\microspace|}}\\
    0_{\phantom{|\microspace\constact\microspace|}} & \id_{|\microspace\constinact\microspace|}
\end{smallmatrix}\right) \Pkj$ \cite[Sec. 6.2d]{themelis_acceleration_2019}, indeed satisfying the idempotence assumption of Theorem \ref{thm:minimization-newton-equiv}.
Substituting this choice of \(\mathcal P_\gamma\) into problem \eqref{eq:minimize-eq-constr}, the constraint can be written as $\dK = -\gamma\, [\Rxnewt]_{\constact.}$
It only remains to solve the, potentially smaller,
unconstrained minimization problem
\newcommand\HJJ{H_{\!\constinact\!\microspace\constinact}}
\begin{equation}\label{eq:minimization-box-unconstr}
    \begin{aligned}[t]
        &\!\microspace\minimize_{\dJ} && \tfrac12 \langle \dJ, \HJJ \dJ \rangle + \langle [ \Rxnewt ]_{\J} + H_{\!\constinact\!\constact} \dK, \dJ \rangle. \\
    \end{aligned}
\end{equation}
Here, \(H_{\!\constinact\!\constact}\) and \(\HJJ\) are the bottom left and right blocks of \(\Pkjt H \Pkj\) respectively.

The structure of \(\Pg\) for the case where the nonsmooth term consists of an \(\ell_1\)-norm is similar to that for box constraints. When $\g = \lambda \Vert \cdot \Vert_1$ with \(\lambda > 0\), we have that $\prox_{\gamma \g}(x)_i = \sign{x_i}\maxbold{(\vert x_i \vert - \gamma \lambda, 0)}$.
Defining $\mathcal{K} \eqdef \left\{ i \in \N_{[1, n]} \mid \vert (\xfb_k)_i \vert \leq \gamma \lambda \right\}$ again yields problem  \eqref{eq:minimization-box-unconstr}.

\subsection{Adaptive regularization}

To deal with cases where \(\HJJ\) is not positive definite, we introduce a radius constraint to keep the solution of \eqref{eq:minimization-box-unconstr} well-defined, giving rise to the following TR subproblem.
\begin{equation}\label{eq:minimization-box-radius-constr}
    \begin{aligned}
        &\minimize_{\dJ}\!\! && \begin{aligned}[t]
            q^\J\!(\dJ) \eqdef\;& \tfrac12 \langle \dJ, \HJJ \dJ \rangle \\
            &+ \langle [\Rxnewt]_\J + H_{\!\constinact\!\constact} \dK, \dJ \rangle
        \end{aligned} \\
        &\stt\!\! && \|\dJ\| \le \trradius \\
    \end{aligned}
\end{equation}
When the radius constraint is inactive, by Theorem \ref{thm:minimization-newton-equiv}, a solution of \eqref{eq:minimization-box-radius-constr}
is exactly a Newton step for the root-finding problem of the residual.
Moreover, when \(\f \in \mathcal{C}^2\), it is also a Newton step for the problem of minimizing the FBE,
as can be verified using \(\nabla\fbe(x) = Q_\gamma(x) R_\gamma(x)\) where \(Q_\gamma = \id - \gamma \nabla^2\f\), and the fact that \(\hat\partial^2\fbe(x) \eqdef Q_\gamma(x) \hat{\partial}R_\gamma(\xnewt)\)
is a LNA for \(\nabla\fbe\) \cite[Cor. 15.4.14]{themelis_acceleration_2019}.

\subsection{Globalization}

To summarize, thus far we described a way to compute update directions for
solving the root-finding problem \(R_\gamma(x^\star) = 0\) by solving a smaller
TR subproblem. Whenever the block \(\HJJ\) of the Hessian is positive
definite and the TR constraint is inactive, that direction is exactly
the Newton direction, and whenever this is not the case, the radius constraint
implicitly regularizes the problem.

Although we can reasonably expect a scheme using those directions to yield fast local convergence, it still lacks a \emph{globalization strategy}.
Hence, we propose to first perform a \emph{forward-backward} (FB) step $\xfb_k \in T_\gamma(x_k)$ at the iterate \(x^k\), and then compute a candidate \emph{accelerated} step $d_k$ as the solution to \eqref{eq:minimization-box-radius-constr} at the point \(\xfb_k\) (i.e., using \(\xnewt = \xfb_k\) when evaluating \(H, \Pg, R_\gamma, \K\) and \(\J\)).
Whether or not the candidate step $d_k$ is accepted is determined by the ratio
\begin{equation} \label{eq:pantr-ratio}
    \trratio_k \eqdef \tfrac{\fbe(\xfb_k) - \fbe(\xfb_k + d_k)}{-q(d_k)},
\end{equation}
which verifies descent on the FBE, similar to \panoc{}. 
It would thus be natural to define a quadratic model of the FBE, e.g.\,using the LNA of \cite[Cor. 15.4.14]{themelis_acceleration_2019} as
\(
    q^{\fbe}\microspace(d) \eqdef \tfrac12 \langle d, \mathcal{Q_\gamma}\lnaRxnewt\,d \rangle + \langle \mathcal{Q_\gamma}R_\gamma(\xnewt), d \rangle
\).
A solution $d_\J^\star$ to the reduced problem \eqref{eq:minimization-box-radius-constr} satisfies $q^\J(d_\J^\star) \leq 0$,
which does not guarantee $q^{\fbe}(d^\star) \leq 0$. 
This is problematic, since a positive $\trratio_k$ need not imply descent on the FBE.
Instead, we define \(
    q(d) \eqdef 
    q^\J\!(\dJ) - \tfrac{1}{2\gamma}\normsq{\dK}
\),
for which we do have that $q(d) \leq q^\J(d_\J)$.
Remark that this model avoids the computation of (expensive) products $H R_\gamma(\xnewt)$.
Based on the ratio \(\trratio_k\), we update the radius $\trradius_k$ using
\begin{equation} \label{eq:radius-update}
    \trradius_{k+1} = \left\{\footnotesize
                \begin{array}{ll}
                    \maxbold\{c_3 \Vert d_k \Vert, \trradius_k\} & \quad \trratio_k \geq \mu_2\\
                    c_2 \trradius_k & \quad \mu_1 \leq \trratio_k < \mu_2\\
                    c_1 \Vert d_k \Vert & \quad \trratio_k < \mu_{1.}
                \end{array}
            \right.
\end{equation}
Since the only purpose of updating $\trradius_k$ is to generate better candidate steps in the next iteration, various alternatives to this heuristic are possible.
When combined, the FB step, the solution of a TR subproblem, and the update of the TR radius make up the \pantr{} scheme, as summarized in Algorithm \ref{alg:pantr}.

\begin{algorithm}
    \small
    \caption{\pantr}
    \label{alg:pantr}
    \begin{algorithmic}
        \Procedure{\pantr$(x_0, \trradius_0, \gamma, \mu_1, \mu_2, c_1, c_2, c_3)$}{}
            \For{$k = 1, 2 \dots$}
                \State Select $\xfb_k \in T_\gamma(x_k)$;
                \State Compute $d_k$ as the solution to \eqref{eq:minimization-box-radius-constr} at \(\xnewt=\xfb_k\);
                \State Compute $\trratio_k$ as in \eqref{eq:pantr-ratio} and update $\trradius_{k+1}$ as in \eqref{eq:radius-update}; \label{alg:fbtr-successful}
                \State $x_{k+1} \leftarrow \left\{
                    \begin{array}{ll}
                        \xfb_k + d_k & \quad \trratio_k \geq \mu_1\\
                        \xfb_k & \quad \trratio_k < \mu_{1;}
                    \end{array}
                \right.$
            \EndFor
        \EndProcedure
    \end{algorithmic}
\end{algorithm}
The reference implementation of \pantr{} solves subproblem \eqref{eq:minimization-box-radius-constr} using the Steihaug CG method.
This method requires only Hessian-vector products and has a limited memory footprint, making it suitable for embedded applications.
Hessian-vector products \(\HJJ\dJ\) can be approximated by a quasi-Newton method on \(\nabla\f\), by finite differences, or using automatic differentiation.
Note that positive definiteness of the quasi-Newton approximants is not required.
The following theorem describes the global subsequential convergence of the presented method.

\begin{theorem}\label{th:pantr-crit}
    Let $\omega(x_k)$ denote the set of cluster points of $\left( x_k \right)_{k \in \N}$.
    Then for the iterates $\left( x_k \right)_{k \in \N}$ of \pantr:
    \begin{enumerate}[label=(\roman*)]
        \item The sequence $\left( \fbe(x_k) \right)_{k \in \N}$ is nonincreasing;
        \item $R_\gamma(x_k) \rightarrow 0$ and $R_\gamma(\xfb_k) \rightarrow 0$ square summably;
        \item $\omega(x_k) = \omega(\xfb_k) \subseteq \fix T_\gamma$;
        \item The sequence $\left( \fbe(x_k) \right)_{k \in \N}$ converges to a finite value $\varphi_\star$, and so does the sequence $\left( \varphi(\xfb_k) \right)_{k \in \N}$ when $\left( x_k \right)_{k \in \N}$ is bounded.
    \end{enumerate}
\end{theorem}
\begin{proof}
    Let $\beta = \frac{1- \gamma \lipschf}{2}$.
    Claim \textit{(i)} follows from the fact that
    \begin{equation} \label{prf:fbe-decrease} \hspace{-0.09cm}
        \!\fbe\microspace(x_{k+1}) \leq \fbe\microspace(\xfb_k) \leq \varphi(\xfb_k) \leq \fbe\microspace(x_k) - \gamma \beta \Vert R_\gamma\microspace(x_k) \Vert^2\!\!,
    \end{equation}
    where we consecutively used \eqref{eq:pantr-ratio} and $q(d_k) \leq 0$, \cite[Prop. 2.2 \textit{(i)}]{stella_forwardbackward_2017}, and \cite[Prop. 2.2 \textit{(ii)}]{stella_forwardbackward_2017}.
    As for \textit{(ii)}, by telescoping \eqref{prf:fbe-decrease} and using the lower-boundedness of $\varphi$ -- and thus of $\fbe$ -- we prove $R_\gamma(x_k) \to 0$ square summably. 
    For $R_\gamma(\xfb_k) \to 0$ square summably, the proof is similar, but instead uses 
    \begin{equation*}
        \fbe(T_\gamma(\xfb_{k})) \leq \varphi(T_\gamma(\xfb_{k})) \leq \fbe(\xfb_{k}) - \gamma \beta \Vert R_\gamma(\xfb_k) \Vert^2\!\!.
    \end{equation*}
    Assume that for some $x' \in \R^n$ and $K \subseteq \N$, we have $\{ x_k \}_{k \in K} \to x'$. 
    Then $\{ \xfb_k \}_{k \in K} \to x'$ as well, since $\Vert \xfb_k - x_k \Vert = \gamma \Vert R_\gamma(x_k) \Vert \to 0$.
    The arbitrarity of $x'$ implies that $\omega(x_k) \subseteq \omega(\xfb_k)$, and by a similar argument also the converse inclusion holds. Hence, $\omega(x_k) = \omega(\xfb_k)$.
    Moreover, $x_k = \prox_{\gamma \g}(x_k - \gamma \nabla \f(x_k)) + \gamma R_\gamma(x_k)$, and since $\{x_k - \gamma \nabla \f(x_k)\}_{k \in K} \to x' - \gamma \nabla \f(x')$, the outer semicontinuity of $\prox_{\gamma \g}$ implies that $x' = \prox_{\gamma \g}(x' - \gamma \nabla \f(x'))$. Hence, $x' \in \fix T_\gamma$, proving \textit{(iii)}.
    From \eqref{prf:fbe-decrease} it follows that $\fbe(x_k) \to \varphi_\star$.
    If $\{ x_k \}_{k \in \N}$ is bounded, then so is $\{ \xfb_k \}_{k \in \N}$ due to compact-valuedness of $\prox_{\gamma \g}$ \cite[Thm. 1.25]{rockafellar_variational_1998}.
    By \cite[Prop. 4.2]{themelis_forward-backward_2018}, $\fbe$ is also $M$-Lipschitz continuous on a compact set containing $\{ x_k \}_{k \in \N}$ and $\{ \xfb_k \}_{k \in \N}$ for some $M > 0$.
    Hence,
    \begin{equation*}
        \begin{aligned}
            \fbe(x_k) - \varphi(\xfb_k) \leq \fbe(x_k) - \fbe(\xfb_k) \leq M \gamma \Vert R_\gamma(x_k) \Vert \to 0.
        \end{aligned}
    \end{equation*}
    Thus we have that $\{ \varphi(\xfb_k) \}_{k \in \N} \to \varphi_\star$, establishing \textit{(iv)}.
\end{proof}

\subsection{Adaptive step size procedure} \label{sec:pantr-adaptive}

To ensure convergence, the step size \(\gamma\) is required to be smaller than
the inverse of the Lipschitz constant $\lipschf$, for example by selecting
\(\gamma = \nicefrac{\alpha}{\lipschf}\) for some \(\alpha \in (0, 1)\).
When the true value of $\lipschf$ is unknown, it can be estimated as follows.
The step size $\gamma$ is updated adaptively by verifying
\begin{equation}\label{eq:Lip-qub-check}
    \f(\xfb_k) \leq \f(x_k)  + \langle \nabla\f(x_k), \xfb_k - x_k \rangle + \tfrac{\alpha}{2 \gamma} \Vert \xfb_k - x_k \Vert^2
\end{equation}
at the start of every iteration of \pantr{}.
Whenever violated, we set $\gamma \leftarrow \nicefrac\gamma2$.
Remark that this can only happen a finite number of times, since when \(\nicefrac \alpha \gamma > \lipschf\), \eqref{eq:Lip-qub-check} is automatically satisfied \cite{themelis_forward-backward_2018}.
Thus, after a finite number of iterations $\gamma$ is constant and all convergence results remain valid as of then.

Moreover, by similar arguments as in \cite{themelis_forward-backward_2018}, \emph{local} Lipschitz-continuity of $\nabla \f$ suffices whenever the candidate directions $d_k$ are bounded and $\g$ has a bounded domain.
The latter is satisfied for $\g = \delta_C$ with $C$ bounded.

\section{Applications and numerical results} \label{sec:numerical}

The effectiveness of our method is demonstrated in the setting of optimal control.
We show that \pantr{} (Alg. \ref{alg:pantr}) performs on par with, or greatly outperforms \alpaqa's existing inner solvers -- \panoc{} \cite{stella_simple_2017,de_marchi_proximal_2022} and structured \panoc{} \cite{pas_alpaqa_2022}.
Run times of \ipopt{} \cite{wachter_implementation_2006}, the default solver in open-source MPC toolboxes such as \texttt{do\_mpc} and \texttt{rockit-meco}, are also reported.

We consider two MPC problems: one with a \emph{hanging chain} model, for which \panoc{} performs well,
and one with a \emph{quadcopter} model and general constraints, selected because it is a problem for which \panoc{} seems to struggle.

\subsection*{Benchmark problems}
Both problems can be stated using
the general formulation
\begin{equation}\tag{OCP} \label{eq:ocp}
    \begin{aligned}
         &  & \minimize_{\mathbf x, \mathbf u} \;\; & \sum_{k=0}^{N-1} \ell_k(x^k, u^k) + \ell_N(x^N)                                                              \\
         &  & \stt\;\;                               & \begin{aligned}[t]
                                                           & x^{k+1} = f(x^k, u^k),                   &  & \forall k \in \N_{[0, N-1]} \\
                                                           & \underline u \le u^k \le \overline u,    &  & \forall k \in \N_{[0, N-1]} \\
                                                           & \underline z \le c(x^k) \le \overline z, &  & \forall k \in \N_{[0, N].}
                                                      \end{aligned}
    \end{aligned}
\end{equation}
\newcommand\nnx{{n_\text{x}}}
\newcommand\nnu{{n_\text{u}}}
\newcommand\nnc{{n_\text{c}}}
The function \(f : \R^\nnx \times \R^\nnu \to \R^\nnx\)
models the discrete-time dynamics of the system.
The matrices \(\mathbf{x} \eqdef (x^1, \dots, x^N) \in \R^{\nnx\times N}\)
and \(\mathbf{u} \eqdef (u^0, \dots, u^{N-1}) \in \R^{\nnu\times N}\) contain the
state and input sequences, respectively. The cost is a sum of the stage costs
\(\ell_k\) and the terminal cost \(\ell_N\). Additionally, the inputs \(u^k\)
are constrained by a rectangular box, and more general state constraints can be
included as well through \(c : \R^\nnx \to \R^\nnc\).
All solvers are applied to the quadcopter and hanging chain OCPs, varying the OCP horizon length from 1 to 60, and simulating the MPC controller for 60 time steps. 
In a first experiment, each solver is cold-started; in a second experiment, 
the solvers are warm-started using the solution and multipliers from the
previous run, shifted by one time step. 

As a first benchmark problem, we consider the \emph{hanging chain model} described in \cite{wirsching_chain_2006}, using the parameters and initial state listed in section III.
The second benchmark problem under consideration is a simplified \emph{quadcopter
model}, governed by the following continuous-time dynamics:
\begin{equation}
    \begin{aligned}
        \dot p = v, \quad\quad \dot v = R(\theta)\,\ttp{(0, 0, a_\text{t})} + g, \quad\quad \dot \theta = \omega.
    \end{aligned}
\end{equation}
The state vector \(x \eqdef (p, v, \theta) \in\R^9\)
consists of the position \(p \in \R^3\), the velocity \(v \in \R^3\),
and the orientation \(\theta \in \R^3\), represented using Euler angles.
The input \(u \eqdef (a_\text{t}, \omega) \in \R^4\) consists of the thrust
\(a_\text{t} \in \R\) and the angular velocity \(\omega \in \R^3\),
and will be determined by the controller as the solution to a finite-horizon
optimal control problem (OCP). \(R(\theta) \in \mathrm{SO}(3)\) represents a rotation matrix, and
\(g = (0, 0, -9.81 \mathrm{m\,s^{-2}})\) is the acceleration due to gravity.
The cost function of the OCP aims to minimize the distance to the target position
\(p_\text{ref} \eqdef (0.25, 0.25, 0.5)\), and penalizes high velocities, angles
and angular velocities, specifically,
\(\ell_k(x, u) \eqdef 10 \normsq{p - p_\text{ref}} + \normsq{v} + \normsq{\theta} + 10 \normsq{\omega} + 10^{-4} a_\text{t}^{2}\) and
\(\ell_N(x) \eqdef 10 \normsq{p - p_\text{ref}} + \normsq{v} + \normsq{\theta}\).
Additionally, we impose constraints on the maximum thrust and on the rate of
rotation by defining \(\underline u \eqdef (0, -0.1, -0.1, -0.1)\) and
\(\overline u \eqdef (49, 0.1, 0.1, 0.1)\). The state
constraints limit the tilt angles and avoid a cylindrical object located at the
origin,\vphantom{\huge M} \(c(x) \eqdef (\theta_x, \theta_y, \cos(\theta_x)\cos(\theta_y), p_x^2 + p_y^2)\),
\(\underline z \eqdef (-\frac\pi2, -\frac\pi2, \cos(\pi/6), 0.1^2)\),
\(\overline z \eqdef (\frac\pi2, \frac\pi2, +\infty, +\infty)\). The system
dynamics are discretized using an explicit fourth-order Runge-Kutta integrator
with a time step \(T_\text{s} = 0.1\,\mathrm{s}\).

\subsection*{Solvers}
A single-shooting formulation of the OCP is used, eliminating the dynamics constraint in \eqref{eq:ocp}.
This yields a NLP of the form \eqref{eq:original-problem}.
Remark that \alpaqa{} thus handles the general state constraint $c$ using its ALM.
However, when no general state constraint is present, \eqref{eq:original-problem} boils down to \eqref{eq:problem-statement} with $\g = \delta_C$, which is solved directly by \panoc{} or \pantr.
We use an L-BFGS buffer of length 50 for the \panoc-based solvers, and for the parameters in Algorithm~\ref{alg:pantr}, we use \(c_1=0.35\), \(c_2=0.99\), \(c_3=10\), \(\mu_1=0.2\), \(\mu_2=0.5\). The step size is determined adaptively as described in Section~\ref{sec:pantr-adaptive}. For \ipopt{}, both the main tolerance and the constraint violation tolerance are set to \(10^{-8}\). The three other solvers declare convergence when \(\norm{\mathbf u-\Pi_{[\underline u, \overline u]}(\mathbf u-\nabla \psi(\mathbf u))}_\infty \le 10^{-8}\)
and \(\norm{c(\mathbf x) - \Pi_{[\underline z, \overline z]}(c(\mathbf x) + \inv \Sigma y)}_\infty \le 10^{-8}\). 
The maximum number of inner iterations per ALM subproblem is set to \(250\), and the initial penalty for the state constraints is set to \(10^4\), with a penalty increase factor of \(5\). The initial inner tolerance is set to \(100\), lowering it by a factor of \(10\) on each ALM iteration.
The necessary derivatives are computed and precompiled using CasADi \cite{andersson_casadi_2019}, except for the hanging chain problem with \ipopt{}, where precompilation is not possible due to the large Hessian matrix, and CasADi's virtual machine (VM) is used for evaluation instead.%
\footnote{The \Cpp{} source code to reproduce the results in this section can be found at
\href{https://github.com/kul-optec/pantr-cdc2023-experiments}{\texttt{github.com/kul-optec/pantr-cdc2023-experiments}}. All experiments were carried out using an Intel Core
i7-11700 CPU at 2.5~GHz.}

\subsection*{Hanging chain results}
\begin{figure}[t]
    \centering
    \includegraphics[width=.95\linewidth, clip, trim=0.4cm 0 0.4cm 0]{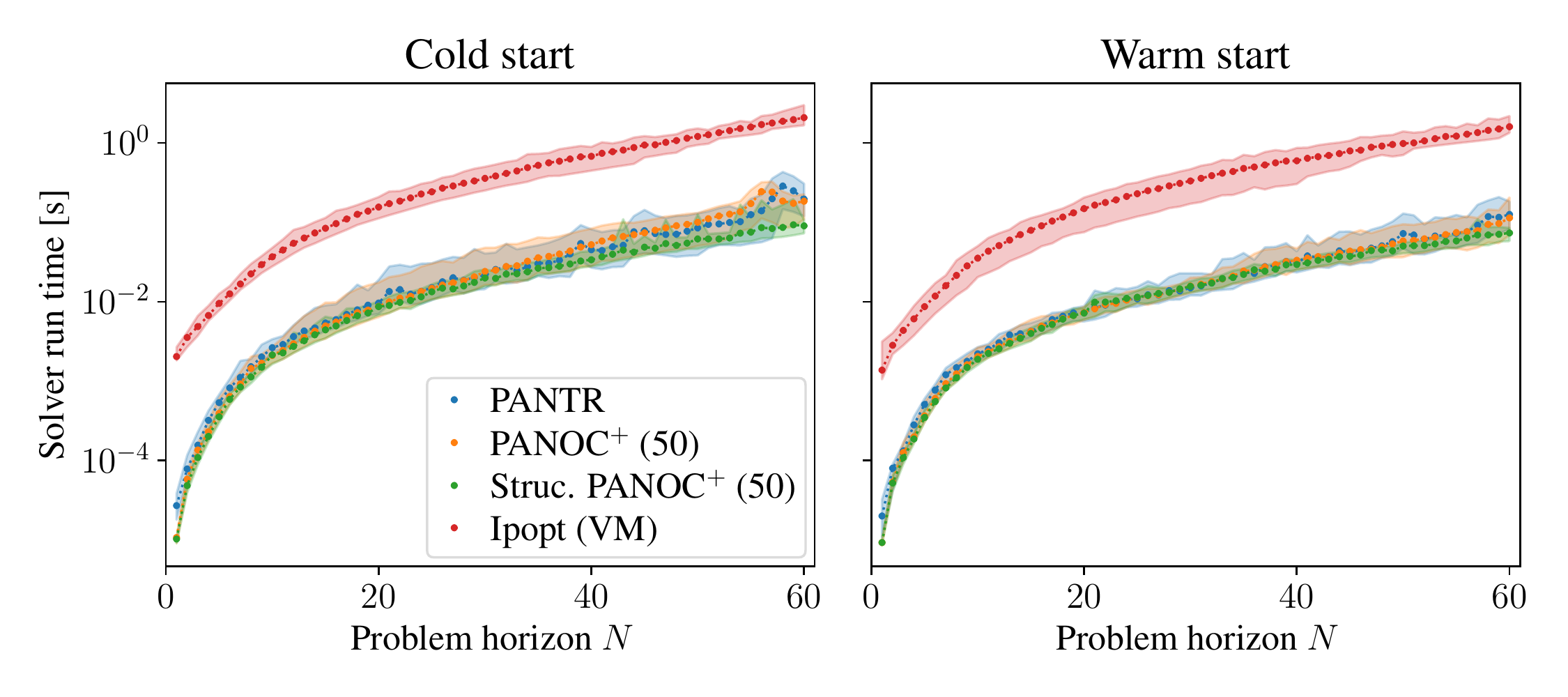}
    \caption{Average solver run times for different horizon lengths of the hanging chain benchmark, with P5/P95 percentiles.}
    \label{fig:simulation-results-chain-avg}
\end{figure}

For the hanging chain problem, where \panoc{} performs quite well, \pantr{} is on par with \panoc{}, and significantly faster than \ipopt{}. 
Fig.~\ref{fig:simulation-results-chain-avg} shows the average run times.
Note that solvers that fully exploit the OCP structure, like the Gauss-Newton variant of \panoc{} from \cite{pas_gauss-newton_2022}, might offer even higher performance.
We do not consider them here, because \pantr{} is applicable to a wider range of problems, whereas \panoc{} with Gauss-Newton only applies to problems with a particular OCP structure.
Similar specialization and exploitation of OCP structure for \pantr{} are the subject of future research.

\begin{figure}[t]
    \centering
    \includegraphics[width=.95\linewidth, clip, trim=0.45cm 0 0.45cm 0]{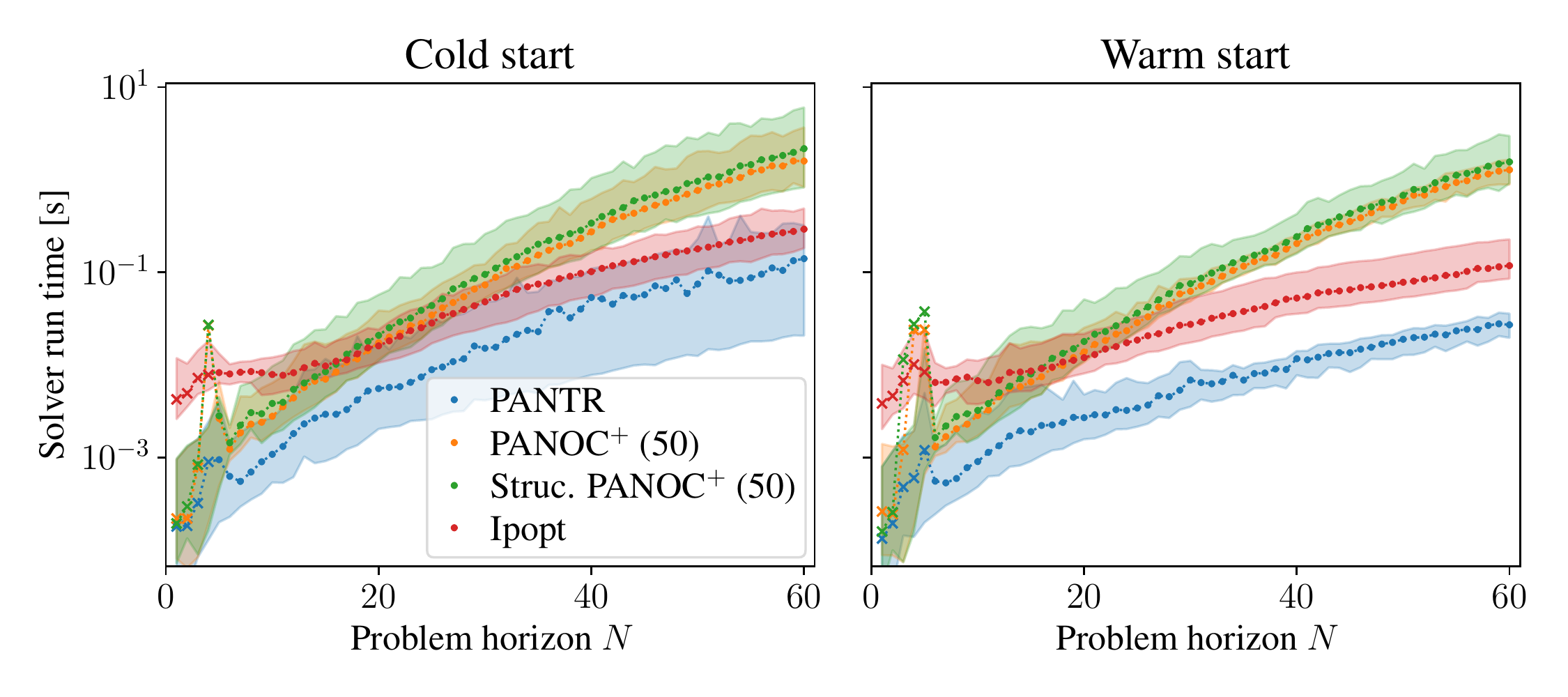}
    \caption{Average solver run times for different horizon lengths of the quadcopter benchmark, with P5/P95 percentiles.}
    \label{fig:simulation-results-quad-avg}
\end{figure}
\begin{figure}[t]
    \centering
    \includegraphics[width=.95\linewidth, clip, trim=0.45cm 0 0.45cm 0]{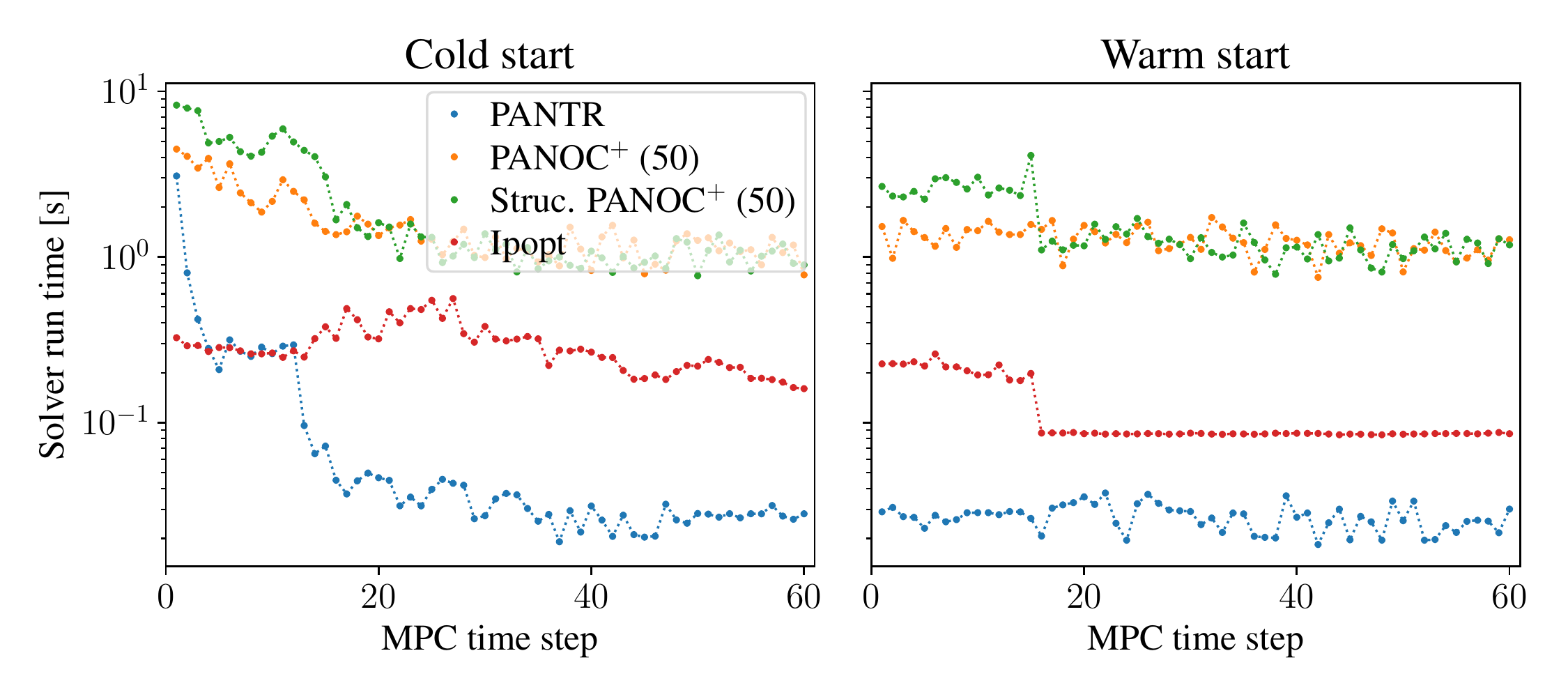}
    \caption{Solver run time per MPC time step, for the quadcopter benchmark with horizon $N=60$.}
    \label{fig:simulation-results-quad-mpc-runtimes}
\end{figure}

\subsection*{Quadcopter results}
Average solver run times for the quadcopter benchmark are shown in Fig.~\ref{fig:simulation-results-quad-avg}. \Panoc{}
and \panoc{} with approximate structured L-BFGS directions do not perform
particularly well here, especially for larger horizons.
This served as one of the motivations for the development of \pantr{}. The latter is able to consistently
outperform \ipopt{}, and scales better with the horizon $N$ than \panoc. Although \ipopt{}
does benefit from warm-starting slightly, the effect is much clearer for \pantr{},
which, when warm-started, achieves average run times that are around three times
faster than \ipopt{}'s for large problems.
\begin{wrapfigure}{r}{0.45\linewidth}%
    \centering
    \includegraphics[width=0.99\linewidth, clip, trim=0.59cm 0.75cm 0.45cm 0.75cm]{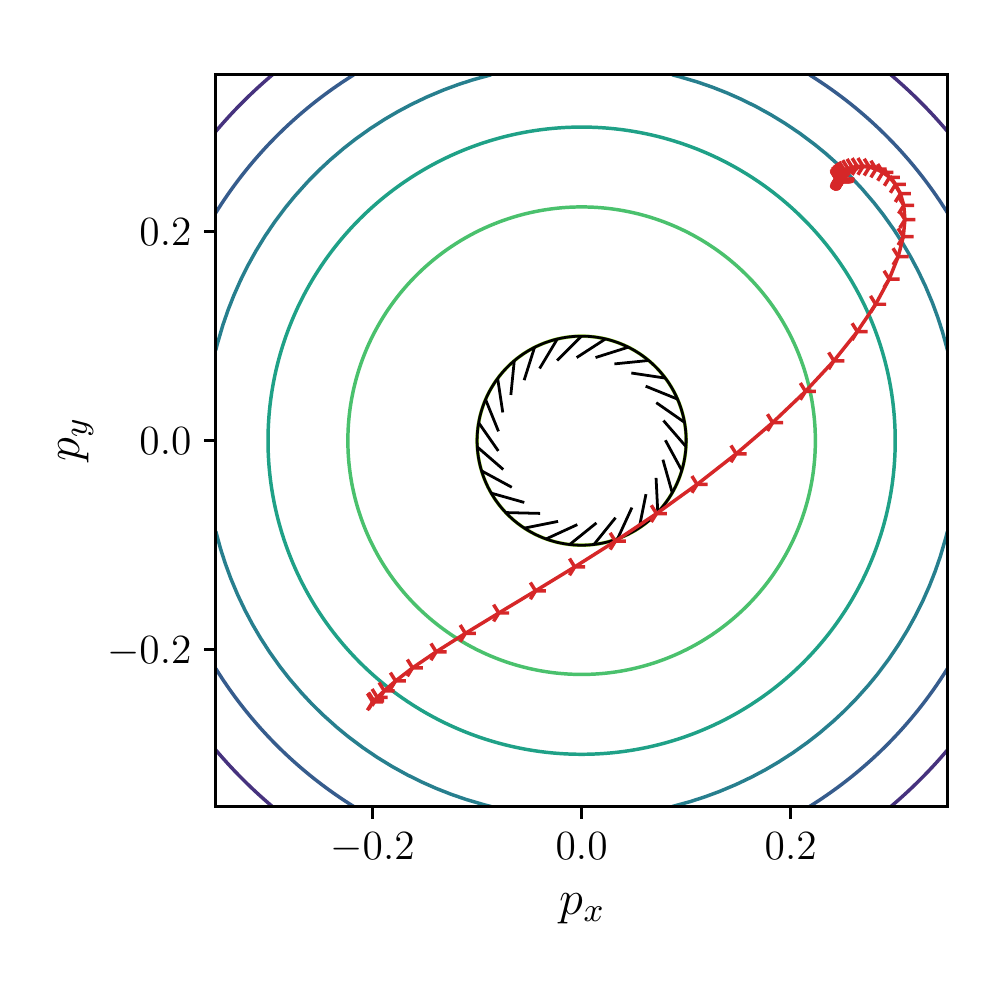}
    \caption{
        Optimal trajectory of quadcopter ($N=60$).
        \vspace{-.5cm}
    }
    \label{fig:quadcopter-solution}
\end{wrapfigure}

Looking at the individual run times for horizon 60 in Fig.~\ref{fig:simulation-results-quad-mpc-runtimes},
it is clear that all solvers take more time during the first 12 MPC
time steps, where the collision constraint is active (as shown in Fig.~\ref{fig:quadcopter-solution}).
For \pantr{}, this is alleviated by warm-starting, resulting in fast convergence throughout.

We observed that the average number of CG iterations per
TR subproblem hovers around 10\% of the number of variables. 
Evaluation of a Hessian-vector product using CasADi is only 1.2 to 3 times more expensive than a
gradient evaluation, further motivating the use of a CG solver for \eqref{eq:minimization-box-radius-constr}.

\section{Conclusion}
This paper presented \pantr, a novel proximal algorithm for nonconvex constrained optimization that is well-suited as an ALM inner solver.
The scheme locally performs Newton steps on an LNA of the fixed-point residual, and as such exploits exact Hessian information of the smooth cost term.
Update directions are computed as solutions to TR subproblems, thus implicitly regularizing the corresponding Newton system.
The presented scheme compares favorably against state-of-the-art NLP solvers in NMPC applications, with the exact second-order information proving particularly beneficial for problems where first-order solvers struggle.

\appendices

\renewcommand*{\bibfont}{\scriptsize}
\printbibliography

\end{document}